\documentclass[11pt]{article} 
\usepackage{amsfonts,amsmath,latexsym,amssymb,mathrsfs,amsthm,comment}
\usepackage{caption}
\usepackage{booktabs}

\let\OLDthebibliography\thebibliography
\renewcommand\thebibliography[1]{
  \OLDthebibliography{#1}
  \setlength{\parskip}{1pt}
  \setlength{\itemsep}{0pt plus 0.0ex}}

\def\numberlikeadb{\global\def\theequation{\thesection.\arabic{equation}}}
\numberlikeadb
\newtheorem{theorem}{Theorem}[section]

\newtheorem{corollary}[theorem]{Corollary}

\newtheorem{proposition}[theorem]{Proposition}
\newtheorem{remark}[theorem]{Remark}

\usepackage{color}

\usepackage{lscape}
\usepackage{caption}
\usepackage{multirow}
\allowdisplaybreaks
\begin{document} 

\title{Inequalities for an integral involving the modified Bessel function of the first kind
}
\author{Robert E. Gaunt\footnote{Department of Mathematics, The University of Manchester, Oxford Road, Manchester M13 9PL, UK}  
} 

\date{} 
\maketitle

\vspace{-10mm}

\begin{abstract} Simple bounds are obtained for the integral $\int_0^x\mathrm{e}^{-\gamma t}t^\nu I_\nu(t)\,\mathrm{d}t$, $x>0$, $\nu>-1/2$, $0\leq\gamma<1$, together with a natural generalisation of this integral. In particular, we obtain an upper bound that holds for all $x>0$, $\nu>-1/2$, $0\leq\gamma<1$, is of the correct asymptotic order as $x\rightarrow0$ and $x\rightarrow\infty$, and possesses a constant factor that is optimal for $\nu\geq0$ and close to optimal for $\nu>-1/2$. We complement this upper bound with several other upper and lower bounds that are tight as $x\rightarrow0$ or as $x\rightarrow\infty$, and apply our results to derive sharper bounds for some expressions that appear in Stein's method for variance-gamma approximation.
\end{abstract}

\noindent{{\bf{Keywords:}}} Modified Bessel function of the first kind; inequality; integral

\noindent{{{\bf{AMS 2010 Subject Classification:}}}} Primary 33C20; 26D15

\section{Introduction}\label{intro}

In the recent papers \cite{gaunt ineq1,gaunt ineq3,gaunt2021}, simple lower and upper bounds, expressed in terms of the modified Bessel function of the first kind $I_\nu(x)$, were derived for the integral
\begin{equation}\label{besint}\int_0^x\mathrm{e}^{-\gamma t}t^\nu I_\nu(t)\,\mathrm{d}t,
\end{equation}
where $x>0$, $\nu>-1/2$, $0\leq\gamma<1$. The conditions on $\nu$ differed from one inequality to the next, but it was always assumed that $\nu>-1/2$ to ensure that the integral exists. For $0 <\gamma < 1$, simple closed-form formulas are not available for this integral. The bounds of \cite{gaunt ineq1,gaunt ineq3,gaunt2021} were needed in the development of Stein's method \cite{chen,np12,stein} for variance-gamma (VG) approximation \cite{aet21,et15,gaunt vg,gaunt vg2,gaunt vg3}. Beyond VG approximation, the bounds
may also be useful in other problems involving modified Bessel functions; for example, inequalities for modified Bessel
functions of the first kind were used by \cite{bs09,baricz3} to establish tight bounds for the generalized Marcum
$Q$-function.

The VG distribution includes as special and limiting cases, amongst others, the classical normal, gamma and Laplace distributions, as well as the product of two standard normal distributions \cite{vg review}.  Stein's method for VG approximation puts existing literature into a more general framework and allows new distributional limits to be considered. Applications include VG approximation for a special case of the $D_2$ statistic from alignment-free sequence comparison, 
 an optimal six moment theorem for VG approximation of double Wiener-It\^{o} integrals, and VG approximation of the generalized Rosenblatt process at extreme critical exponent. 
 Related recent works
 in which the VG distribution arises as a limit distribution include \cite{aaps17,ag18,azmooden,bt17}.

At the heart of Stein's method for VG approximation lies the function $f_h:\mathbb{R}\rightarrow\mathbb{R}$ defined by
\begin{align}
\label{vgsoln}
f_h(x)&=-\frac{\mathrm{e}^{-\beta x} K_{\nu}(|x|)}{|x|^{\nu}} \int_0^x \mathrm{e}^{\beta t} |t|^{\nu} I_{\nu}(|t|) h(t) \,\mathrm{d}t  -\frac{\mathrm{e}^{-\beta x} I_{\nu}(|x|)}{|x|^{\nu}} \int_x^{\infty} \mathrm{e}^{\beta t} |t|^{\nu} K_{\nu}(|t|)h(t)\,\mathrm{d}t, 
\end{align}
where $x\in\mathbb{R}$, $\nu>-1/2$, $-1<\beta<1$, and $h:\mathbb{R}\rightarrow\mathbb{R}$ has mean zero with respect to the  VG probability measure. Here, $K_\nu(x)$ is the modified Bessel function of second kind. The function $f_h$ is the solution of the so-called \emph{Stein equation}, and in order to implement Stein's method for VG approximation, uniform bounds, in terms of supremum norms of $h$ and its lower order derivatives, are required on the solution $f_h$ and its first four derivatives. To achieve such bounds, uniform bounds were required for a number of expressions involving the modified Bessel functions $I_\nu(x)$ and $K_\nu(x)$. Most of these expressions were bounded in \cite{gaunt ineq1,gaunt ineq2, gaunt ineq3}; however, establishing uniform bounds on the quantities
\begin{equation}\label{terms} M_{\nu,\beta,n}(x)=\frac{\mathrm{e}^{-\beta x}K_{\nu+n}(x)}{x^{\nu-1}}\int_0^x \mathrm{e}^{\beta t}t^{\nu}I_\nu(t)\,\mathrm{d}t, \quad n=0,1,2,
\end{equation}
that hold for all $\nu>-1/2$, $-1<\beta<1$, $x\geq0$ (it suffices to consider $x\geq0$) proved more troublesome. In a recent paper \cite{gaunt2021}, uniform bounds were obtained for the expressions in (\ref{terms}), which allowed for a technical advance in Stein's method for VG approximation. 

As a consequence of the uniform bound for $M_{\nu,\beta,2}(x)$, \cite{gaunt2021} was able to solve an open problem of \cite{gaunt ineq3}, which was to find a constant $C_{\nu,\gamma}>0$ such that, for all $x>0$,
\begin{equation}\label{target} \int_0^x \mathrm{e}^{-\gamma t}t^\nu I_\nu(t)\,\mathrm{d}t<C_{\nu,\gamma}\mathrm{e}^{-\gamma x}x^\nu I_{\nu+1}(x), \quad \nu>-\tfrac{1}{2},\:0\leq\gamma<1.
\end{equation}
In doing so, \cite[inequality (2.12)]{gaunt2021} obtained the first bound in the literature for the integral (\ref{besint}) that is valid for all $x>0$, $\nu>-1/2$, $0\leq\gamma<1$ and is of the correct asymptotic order as $x\rightarrow0$ and $x\rightarrow\infty$. However, the indirect approach to bounding the integral (\ref{besint}) resulted in the 
rather unsatisfactory 
estimate $C_{\nu,\gamma}\leq (4\nu+14)/((2\nu+1)(1-\gamma))$. Indeed, as a means of comparison, inequality (2.19) of \cite{gaunt ineq3} states that, if we restrict to $\nu\geq1/2$ and $0\leq\gamma<1$, then, for $x>0$,
\begin{align}\int_0^x \mathrm{e}^{-\gamma t}t^{\nu}I_\nu(t)\,\mathrm{d}t&<\frac{\mathrm{e}^{-\gamma x}x^\nu}{(2\nu+1)(1-\gamma)}\Big(2(\nu+1)I_{\nu+1}(x)-I_{\nu+3}(x)\Big)\label{baaad}\\
\label{gau1}&<\frac{2(\nu+1)}{(2\nu+1)(1-\gamma)}\mathrm{e}^{-\gamma x}x^\nu I_{\nu+1}(x),
\end{align}
with the same inequalities being valid for all $\nu>-1/2$ in the case $\gamma=0$ (see inequalities (2.17) and (2.18) of \cite{gaunt ineq3}). 



In this paper, we take a more direct approach to bounding the integral (\ref{besint}), which allows us to obtain a bound of the form (\ref{target}) that holds for all $x>0$, $\nu>-1/2$, $0\leq\gamma<1$ with a constant factor $C_{\nu,\gamma}$ that significantly improves on that of \cite[inequality (2.12)]{gaunt2021} (see Theorem \ref{thm1}). 
As mentioned in Remark \ref{rem1}, for $\nu\geq0$ the constant given in the bound of Theorem \ref{thm1} is the same as in inequality (\ref{gau1}), and we therefore extend the range of validity of the inequality from $\nu\geq1/2$ to $\nu\geq0$. 

In Theorems \ref{thm2} and \ref{cor1}, we complement our upper bound from Theorem \ref{thm1} with several upper and lower bounds for integrals of a similar form to (\ref{besint}), one of which is required in the proof of Theorem \ref{thm1} and some of which are derived using our upper bound from Theorem \ref{thm1}. In particular, in Theorem \ref{thm2}, we derive sharp upper and lower bounds for a natural generalisation of the integral (\ref{besint}).
We also apply our upper bound from Theorem \ref{thm1} to deduce uniform bounds for the expressions in (\ref{terms}) (see Corollary \ref{cor2}). It turns out that our approach of first bounding the integral in fact also leads to better constants here, which will allow for improved constants to be given in a number of the bounds for the solution (\ref{vgsoln}) of the VG Stein equation and for derivatives of the solution. 


\section{Main results}\label{sec2}


The following theorem is the main result of this paper. 

\begin{theorem}\label{thm1}Let $0\leq\gamma<1$ and $\nu>-1/2$. Then, for $x>0$,
\begin{align}\label{mainbd}\int_0^x \mathrm{e}^{-\gamma t}t^{\nu}I_{\nu}(t)\,\mathrm{d}t<\frac{2(\nu+1)+c_\nu}{(2\nu+1)(1-\gamma)}\mathrm{e}^{-\gamma x}x^\nu I_{\nu+1}(x),
\end{align}
where $c_\nu=\max\{0,\,-4\nu(\nu+1)\}$.

\end{theorem}

\begin{remark}\label{rem1}Since $c_\nu=0$ for $\nu\geq0$, inequality (\ref{mainbd}) extends the range of validity of inequality (\ref{gau1}) from $\nu\geq1/2$ to $\nu\geq0$. Moreover, since $c_\nu<1$ for $\nu>-1/2$, for the full range of parameters $x>0$, $\nu>-1/2$, $0\leq\gamma<1$ we have the simple bound
\begin{equation*}
\int_0^x \mathrm{e}^{-\gamma t}t^{\nu}I_{\nu}(t)\,\mathrm{d}t<\frac{2\nu+3}{(2\nu+1)(1-\gamma)}\mathrm{e}^{-\gamma x}x^\nu I_{\nu+1}(x),    
\end{equation*}
in which the estimate for the constant $C_{\nu,\gamma}$ significantly improves on that of \cite[inequality (2.12)]{gaunt2021}. 
\end{remark}


In the following theorem, we provide upper and lower bounds for a natural generalisation of the integral (\ref{besint}). These bounds complement inequalities given in \cite[Theorem 2.1]{gaunt ineq1}  and \cite[Theorem 2.3]{gaunt ineq3} for this integral, as well as the bounds of \cite[Theorem 2.6]{gaunt ineq6} for the integral $\int_0^x \mathrm{e}^{-\gamma t}t^{-\nu} I_{\nu+n}(t)\,\mathrm{d}t$.

\begin{theorem}\label{thm2} Let $0\leq\gamma<1$ and suppose that $\nu,n\in\mathbb{R}$ are such that $\nu>-(n+1)/2$. 

\vspace{2mm}

\noindent  (i) Suppose further that $n>-1$ and $\nu\geq1/2$. Then, for $x>0$,
\begin{align} \label{new1}
\int_0^x \mathrm{e}^{-\gamma t}t^{\nu}I_{\nu+n}(t)\,\mathrm{d}t&<\frac{\mathrm{e}^{-\gamma x}x^\nu}{(2\nu+n+1)(1-\gamma)}\Big(2(\nu+n+1)I_{\nu+n+1}(x) -(n+1)I_{\nu+n+3}(x)\Big). \end{align}
When $\gamma=0$ inequality (\ref{new1}) holds without the need for the restriction $\nu\geq1/2$ (that is it holds provided $n>-1$ and $\nu>-(n+1)/2$).
Inequality (\ref{new1}) becomes an equality when $\gamma=0$ and $n=-1$, $\nu>0$, and the inequality is reversed when $\gamma=0$ and $-3<n<-1$, $\nu>-(n+1)$ (again without the need for the additional condition that $\nu\geq1/2$).

\vspace{2mm}

\noindent (ii) Suppose that $n>-1$ and $\nu>-(n+1)/2$. Then, for $x>0$,
\begin{align}\label{lower4} \int_0^x \mathrm{e}^{-\gamma t}t^\nu I_{\nu+n}(t)\,\mathrm{d}t&>\frac{\mathrm{e}^{-\gamma x}x^\nu}{2\nu+n+1}\Big(2(\nu+n+1)I_{\nu+n+1}(x)\nonumber\\
&\quad -\frac{2(n+1)(\nu+n+3)}{2\nu+n+3}I_{\nu+n+3}(x)+\frac{(n+1)(n+3)}{2\nu+n+3}I_{\nu+n+5}(x)\Big). 
\end{align}
\noindent (iii) Inequality (\ref{new1}) is tight as $x\rightarrow\infty$, and is also tight as $x\rightarrow0$ when $\gamma=0$. Inequality (\ref{lower4}) is tight as $x\rightarrow0$, and is also tight as $x\rightarrow\infty$ when $\gamma=0$.
\end{theorem}

In the following theorem, we state some further lower bounds for integrals of a similar form to the integral (\ref{besint}). One of these lower bounds will be used in the proof of Theorem \ref{thm1}.

\begin{theorem}\label{cor1} Let $0\leq\gamma<1$. Then, for $x>0$,
\begin{align}\label{lower1} \int_0^x \mathrm{e}^{-\gamma t}t^{\nu+1}I_{\nu}(t)\,\mathrm{d}t&>\mathrm{e}^{-\gamma x}x^{\nu+1}\sum_{k=0}^\infty\gamma^kI_{\nu+k+1}(x), \quad \nu>-1, \\
\label{intineq0}\int_0^x \mathrm{e}^{-\gamma t}t^{\nu}I_{\nu+1}(t)\,\mathrm{d}t&>\frac{1}{1-\gamma}\bigg\{1-\frac{2\nu(2\nu+c_{\nu-1})}{(2\nu-1)(1-\gamma)}\frac{1}{x}\bigg\}\mathrm{e}^{-\gamma x}x^\nu I_\nu(x), \quad \nu>\tfrac{1}{2}.
\end{align}
From inequality (\ref{intineq0}) we deduce the lower bound
\begin{equation}\label{lower2}
\int_0^x \mathrm{e}^{-\gamma t}t^{\nu}I_{\nu}(t)\,\mathrm{d}t>\frac{1}{1-\gamma}\bigg\{1-\frac{2\nu(2\nu+c_{\nu-1})}{(2\nu-1)(1-\gamma)}\frac{1}{x}\bigg\}\mathrm{e}^{-\gamma x}x^\nu I_\nu(x), \quad \nu>\tfrac{1}{2}.    
\end{equation}
Inequalities (\ref{lower1})--(\ref{lower2}) are tight as $x\rightarrow\infty$. 
\end{theorem}

Theorem \ref{thm1} extends the range of validity of inequality (\ref{gau1}) from $\nu\geq1/2$ to $\nu\geq0$.
In the same spirit, it is natural to ask whether the range of validity of the following simple bound of \cite[inequality (2.5)]{gaunt ineq1} can extended: Let $0\leq\gamma<1$ and $\nu\geq1/2$. Then, for $x>0$,
\begin{equation}
\int_0^x \mathrm{e}^{-\gamma t}t^{\nu}I_{\nu}(t)\,\mathrm{d}t<\frac{1}{1-\gamma}\mathrm{e}^{-\gamma x}x^\nu I_\nu(x).  \label{gau2}  
\end{equation}
(Inequality (\ref{gau2}) is stated in \cite{gaunt ineq1} as a weak inequality, but the inequality is actually strict.) One may suspect that the range of validity of inequality (\ref{gau2}) can be extended from $\nu\geq1/2$ to $\nu>-1/2$. However, it turns out that the range of validity $\nu\ge1/2$ is in fact optimal, as can be seen as a special case of the following proposition.

\begin{proposition}\label{prop1} Suppose that $0\leq\gamma<1$, and that $\mu,\nu\in\mathbb{R}$ are such that $\mu+\nu>-1$. 

\vspace{2mm}

\noindent (i) Suppose further that $\mu\geq\nu\geq1/2$. Then, for $x>0$,
\begin{equation}
\int_0^x \mathrm{e}^{-\gamma t}t^{\mu}I_{\nu}(t)\,\mathrm{d}t<\frac{1}{1-\gamma}\mathrm{e}^{-\gamma x}x^\mu I_\nu(x).  \label{gau222}
\end{equation}
\noindent (ii) If $\mu<1/2$, then there exists $x*>0$ such that for all $x\geq x_*$ we have
\begin{equation}\label{lar}
\int_0^x \mathrm{e}^{-\gamma t}t^{\mu}I_{\nu}(t)\,\mathrm{d}t>\frac{1}{1-\gamma}\mathrm{e}^{-\gamma x}x^\mu I_\nu(x). 
\end{equation}
    
\end{proposition}

In the following corollary, we obtain sharper bounds for the expressions in (\ref{terms}) than are given in Theorem 3.2 of \cite{gaunt2021}. These bounds are derived using our upper bound (\ref{mainbd}) for the integral (\ref{besint}).

\begin{corollary}\label{cor2}Suppose that $-1<\beta\leq0$ and $\nu>-1/2$. Then, for $x\geq0$,
\begin{align}\label{c1}\frac{\mathrm{e}^{-\beta x}K_{\nu+2}(x)}{x^{\nu-1}}\int_0^x \mathrm{e}^{\beta t}t^{\nu}I_\nu(t)\,\mathrm{d}t&<\frac{2(\nu+1)+c_\nu}{(2\nu+1)(1+\beta)},\\
\label{c2}\frac{\mathrm{e}^{-\beta x}K_{\nu+1}(x)}{x^{\nu-1}}\int_0^x \mathrm{e}^{\beta t}t^{\nu}I_\nu(t)\,\mathrm{d}t&<\frac{\nu+1+c_\nu/2}{(2\nu+1)(1+\beta)}, \\
\label{c3}\frac{\mathrm{e}^{-\beta x}K_{\nu}(x)}{x^{\nu-1}}\int_0^x \mathrm{e}^{\beta t}t^{\nu}I_\nu(t)\,\mathrm{d}t&<\frac{\nu+1+c_\nu/2}{(2\nu+1)(1+\beta)}.
\end{align}
\end{corollary}

\begin{remark} 1. Inequality (\ref{lower1}) is used in the proof of Theorem \ref{thm1}. It complements the following lower bound of \cite[inequality (2.10)]{gaunt2021} for the integral (\ref{besint}), which is also used in the proof of Theorem \ref{thm1}: Let $0\leq\gamma<1$ and $\nu>-1/2$. Then, for $x>0$,
\begin{equation}
\label{lower3}  \int_0^x \mathrm{e}^{-\gamma t}t^{\nu}I_{\nu}(t)\,\mathrm{d}t>\mathrm{e}^{-\gamma x}x^{\nu}\sum_{k=0}^\infty\gamma^kI_{\nu+k+1}(x).  
\end{equation}
We also note that inequality (\ref{lower1}) was established in the case $\gamma=0$ by \cite[inequality (2.20)]{gaunt ineq3}. Inequalities (\ref{intineq0}) and (\ref{lower2}) are deduced from our upper bound (\ref{mainbd}). Inequality (\ref{lower2}) extends the range of validity of inequality (2.9) of \cite{gaunt2021} from $\nu\geq3/2$ to $\nu\geq1$ (since $c_{\nu-1}=0$ for $\nu\geq1$, in which case both bounds are equal), whilst inequality (\ref{intineq0}) is a strict improvement, since $I_{\nu+1}(x)<I_\nu(x)$ for $x>0$, $\nu>1/2$ (see inequality (\ref{Imon})). 

\vspace{2mm}

\noindent 2. Inequality (\ref{new1}) is a natural generalisation of inequality (\ref{gau1}) of \cite{gaunt ineq3}. Inequality (\ref{gau1}) was derived using the classical inequality $I_{\nu+1}(x)<I_{\nu}(x)$,  $x>0$, $\nu\geq1/2$ (see (\ref{Imon})), whilst inequality (\ref{new1}) is derived using a recent interesting refinement of this inequality (see inequality (\ref{ytineq})). Inequality (\ref{new1}) can also be viewed as a generalisation of inequality (2.17) of \cite{gaunt ineq3}, which was derived for the case $\gamma=0$, $\nu>-1/2$. 

When restricted to the case $\gamma=0$, the lower bound (\ref{lower4}) is the first lower bound in the literature for the integral $\int_0^x t^\nu I_{\nu}(t)\,\mathrm{d}t$ and more generally for  the integral $\int_0^x t^\nu I_{\nu+n}(t)\,\mathrm{d}t$ that is tight as $x\rightarrow0$ and as $x\rightarrow\infty$. This inequality complements the upper bound (\ref{new1}) with $\gamma=0$, which is also tight as $x\rightarrow0$ and as $x\rightarrow\infty$. Taken together these bounds yield the following two-sided inequality that is tight  as $x\rightarrow0$ and as $x\rightarrow\infty$: Let $n>-1$ and $\nu>-(n+1)/2$. Then, for $x>0$, 
\begin{equation}\label{twos}
L_{\nu,n}(x)<\int_0^x t^\nu I_{\nu+n}(t)\,\mathrm{d}t<U_{\nu,n}(x),    
\end{equation}    
where
\begin{align*}
L_{\nu,n}(x)&=\frac{x^\nu}{2\nu+n+1}\Big(2(\nu+n+1)I_{\nu+n+1}(x)\nonumber\\
&\quad -\frac{2(n+1)(\nu+n+3)}{2\nu+n+3}I_{\nu+n+3}(x)+\frac{(n+1)(n+3)}{2\nu+n+3}I_{\nu+n+5}(x)\Big),\\
U_{\nu,n}(x)&=\frac{x^\nu}{2\nu+n+1}\Big(2(\nu+n+1)I_{\nu+n+1}(x)-(n+1)I_{\nu+n+3}(x)\Big).
\end{align*}

We used \emph{Mathematica} to calculate the relative error in estimating $F_{\nu}(x)=\int_0^x t^\nu I_{\nu}(t)\,\mathrm{d}t$ by the lower and upper bounds $L_{\nu,0}(x)$ and $U_{\nu,0}(x)$  for a range of values of $x$ and $\nu$. The results are reported in Tables \ref{table1} and \ref{table2}. As one might expect given that the two-sided inequality (\ref{twos}) is tight as $x\rightarrow0$ and as $x\rightarrow\infty$, we observe that, for fixed $\nu$, the relative error in approximating the integral $F_\nu(x)$ by either $L_{\nu,0}(x)$ or $U_{\nu,0}(x)$ initially increases (from rather small values when $x=1$) before decreasing as $x$ get larger. We also see that, for fixed $x$, the relative error in approximating the integral $F_\nu(x)$ by either $L_{\nu,0}(x)$ or $U_{\nu,0}(x)$ decreases as $\nu$ increases.

\begin{table}[h]
  \centering
  \caption{Relative error in approximating the integral $F_\nu(x)$ by the lower bound  $L_{\nu,0}(x)$.
  }
\label{table1}
\normalsize{
\begin{tabular}{l*{8}{c}}
\hline
& \multicolumn{8}{c}{$x$} \\
\cmidrule(lr){2-9}
$\nu$ & 1 &2.5 & 5 & 10 & 15 & 25 & 50 & 100 \\
\hline
        $-0.25$ &0.0006  &0.0199  &0.1528  &0.3593  &0.3747   &0.3105  &0.1943 &0.1081 \\ 
        0 & 0.0002 & 0.0074 &0.0528  &0.1305  &0.1425   &0.1227 &0.0789 &0.0445  \\ 
        1 & 0.0000 & 0.0006 & 0.0046 &0.0154  &0.0199   &0.0199  &0.0142 &0.0085 \\ 
        2.5 & 0.0000 & 0.0000 &0.0005  &0.0023  &0.0037   &0.0045  &0.0038 & 0.0025 \\
        5 &0.0000 & 0.0000 & 0.0000 & 0.0003 &0.0006 &0.0009   &0.0010 &0.0007 
        \\ \hline
    \end{tabular}}
\end{table} 

\begin{table}[h]
  \centering
  \caption{Relative error in approximating the integral $F_\nu(x)$ by the upper bound  $U_{\nu,0}(x)$.
  }
\label{table2}
\normalsize{
\begin{tabular}{l*{8}{c}}
\hline
& \multicolumn{8}{c}{$x$} \\
\cmidrule(lr){2-9}
$\nu$ & 1 & 2.5 & 5 & 10 & 15 & 25 & 50 & 100 \\
\hline
        $-0.25$ &0.0403  &0.2132  &0.4675  & 0.4323  &0.3268   & 0.2137 & 0.1134 & 0.0584 \\ 
        0 & 0.0199 & 0.0991 & 0.2038 &  0.1973 & 0.1543  & 0.1034 & 0.0558 & 0.0290  \\ 
        1 & 0.0030 & 0.0156 & 0.0368 &0.0464  &0.0411   &0.0303  &0.0175 & 0.0094 \\ 
        2.5 & 0.0001 &0.0030 & 0.0084 &0.0144  &0.0149  &0.0125   &0.0080  &0.0045 \\
        5 &0.0000 &0.0005 & 0.0017 &0.0039  & 0.0049 &0.0050   &0.0037  &0.0023 
        \\ \hline
    \end{tabular}}
\end{table} 

\end{remark}

\section{Proofs of main results}\label{sec3}

We first prove Theorem \ref{cor1}, as inequality (\ref{lower1}) is required in the proof of Theorem \ref{thm1}. 

\vspace{3mm}

\noindent{\emph{Proof of Theorem \ref{cor1}.}}   By an integration by parts using the differentiation formula (\ref{diffone}) we obtain that
\begin{align}\label{jj25}\int_0^x\mathrm{e}^{-\gamma t}t^{\nu+1}I_\nu(t)\,\mathrm{d}t=\mathrm{e}^{-\gamma x}x^{\nu+1}I_{\nu+1}(x)+\gamma\int_0^x \mathrm{e}^{-\gamma t}t^{\nu+1}I_{\nu+1}(t)\,\mathrm{d}t,
\end{align}
where we made use of the limit $\lim_{x\rightarrow0}x^{\nu+1}I_{\nu+1}(x)=0$ for $\nu>-1$ (see (\ref{Itend0})). On using the integral inequality  (\ref{lower3}) we now get that
\begin{align*}
\int_0^x\mathrm{e}^{-\gamma t}t^{\nu+1}I_\nu(t)\,\mathrm{d}t&>\mathrm{e}^{-\gamma x}x^{\nu+1}I_{\nu+1}(x)+\gamma\mathrm{e}^{-\gamma x}x^{\nu+1}\sum_{j=0}^\infty\gamma^j I_{\nu+j+2}(x)\\
&=\mathrm{e}^{-\gamma x}x^{\nu+1}\sum_{k=0}^\infty\gamma^k I_{\nu+k+1}(x).
\end{align*}


\noindent (ii) Replacing $\nu$ by $\nu-1$ in equation (\ref{jj25}), using identity (\ref{Iidentity}) and then rearranging gives that
\begin{align}\label{1stint}\int_0^x \mathrm{e}^{-\gamma t}t^{\nu}I_{\nu+1}(t)\,\mathrm{d}t-\gamma \int_0^x \mathrm{e}^{-\gamma t}t^{\nu}I_{\nu}(t)\,\mathrm{d}t=\mathrm{e}^{-\gamma x}x^\nu I_\nu(x)-2\nu\int_0^x \mathrm{e}^{-\gamma t}t^{\nu-1}I_{\nu}(t)\,\mathrm{d}t,
\end{align}
with the equality valid for $\nu>0$ since $\lim_{x\rightarrow0}x^{\nu}I_{\nu}(x)=0$ for $\nu>0$. Applying inequality (\ref{Imon}) to the second and third  integrals in (\ref{1stint}) and rearranging yields
\begin{align}\int_0^x \mathrm{e}^{-\gamma t}t^{\nu}I_{\nu+1}(t)\,\mathrm{d}t>\frac{1}{1-\gamma}\bigg\{\mathrm{e}^{-\gamma x}x^\nu I_\nu(x)-2\nu\int_0^x \mathrm{e}^{-\gamma t}t^{\nu-1}I_{\nu-1}(t)\,\mathrm{d}t\bigg\}, \label{1stint0}
\end{align}
where the second integral in (\ref{1stint0}) exists for $\nu>1/2$ (see (\ref{Itend0})). Inequality (\ref{intineq0}) now follows on applying the bound (\ref{mainbd}) to inequality (\ref{1stint0}). 

\vspace{2mm}

\noindent (iii) Apply inequality (\ref{Imon}) to the integral in inequality (\ref{intineq0}). 

\vspace{2mm}

\noindent (iv) By a simple asymptotic analysis using the limiting form (\ref{Itendinfinity}) we have that, for $\nu>-1/2$ and $0\leq\gamma<1$, 
\begin{equation}\label{eqeq1} \int_0^x \mathrm{e}^{-\gamma t}t^\nu  I_{\nu}(t)\,\mathrm{d}t\sim \frac{1}{\sqrt{2\pi}(1-\gamma)}x^{\nu-1/2}\mathrm{e}^{(1-\gamma)x}, \quad x\rightarrow\infty,
\end{equation}
and we also have from (\ref{Itendinfinity}) that, for $n\in\mathbb{R}$,
\begin{equation}\label{eqeq2}\mathrm{e}^{-\gamma x}x^\nu I_{\nu+n}(x)\sim  \frac{1}{\sqrt{2\pi}}x^{\nu-1/2}\mathrm{e}^{(1-\gamma)x}, \quad x\rightarrow\infty.
\end{equation}
That inequalities (\ref{lower1})--(\ref{lower2}) are tight as $x\rightarrow\infty$ is immediate from the limiting forms (\ref{eqeq1}) and (\ref{eqeq2}). In the case of the verification that inequality (\ref{lower1}) is tight as $x\rightarrow1$ we make use of the geometric series formula $\sum_{k=0}^\infty\gamma^k=1/(1-\gamma)$ for $0\leq\gamma<1$.
\qed

\vspace{3mm}

We are now ready to prove Theorem \ref{thm1}. The proof
is split into three parts. It is interesting to note  that the analysis in the first part of the proof is sufficient to prove the inequality in the case $\nu\geq1/2$, $0\leq\gamma<1$. We therefore provide an alternative derivation of inequality (2.19) of \cite{gaunt ineq3}, which established inequality (\ref{mainbd}) in this restricted case.

\vspace{3mm}

\noindent{\emph{Proof of Theorem \ref{thm1}.}} 1. Consider the function
\begin{equation*}
u_{\nu,\gamma}(x)=\frac{2(\nu+1)}{(2\nu+1)(1-\gamma)}\mathrm{e}^{-\gamma x}x^\nu I_{\nu+1}(x)-\int_0^x \mathrm{e}^{-\gamma t}t^{\nu}I_{\nu}(t)\,\mathrm{d}t.   
\end{equation*}
We will begin by proving that, for fixed $\nu>-1/2$ and $0\leq\gamma<1$, $u_{\nu,\gamma}(x)$ is an increasing function of $x$ for all $x> x_*$ for some $x_*>0$, which we will later identify. This will reduce the problem of proving inequality (\ref{mainbd}) for all $x>0$ to proving the inequality for all $x\in(0, x_*]$. As we will see, this analysis will in fact prove inequality (\ref{mainbd}) for all $x>0$ in the case $\nu\geq1/2$, $0\leq\gamma<1$.

By applying the differentiation formula (\ref{diffone}) we obtain that
\begin{align*}u_{\nu,\gamma}'(x)&=\frac{2(\nu+1)}{(2\nu+1)(1-\gamma)}\frac{\mathrm{d}}{\mathrm{d}x}\big(\mathrm{e}^{-\gamma x}x^{-1}\cdot x^{\nu+1} I_{\nu+1}(x)\big)- \mathrm{e}^{-\gamma x}x^\nu I_\nu(x)\\
&=\frac{2(\nu+1)}{(2\nu+1)(1-\gamma)}\mathrm{e}^{-\gamma x}x^\nu\big(I_{\nu}(x)-x^{-1}I_{\nu+1}(x)-\gamma I_{\nu+1}(x)\big)-\mathrm{e}^{-\gamma x}x^\nu I_\nu(x) \\
&=\frac{2(\nu+1)}{(2\nu+1)(1-\gamma)}\mathrm{e}^{-\gamma x}x^\nu I_\nu(x)\bigg(\frac{1+\gamma(2\nu+1)}{2(\nu+1)}-(x^{-1}+\gamma)\frac{I_{\nu+1}(x)}{I_\nu(x)}\bigg).
\end{align*}
To find an $x_*$ such that $u_{\nu,\gamma}'(x)>0$ for all $x>x_*$ it therefore suffices to find an $x_*$ such that for all $x>x_*$ the following inequality holds
\begin{align}\label{cond1}
\frac{I_{\nu+1}(x)}{xI_\nu(x)}< \frac{1+\gamma(2\nu+1)}{2(\nu+1)(1+\gamma x)}.   
\end{align}
We now recall the following simple bound of \cite{nasell2}: for $x>0$ and $\nu>-1/2$,
\begin{align}
\frac{I_{\nu+1}(x)}{I_\nu(x)}&<\frac{x}{\nu+1/2+x}. \label{nas2}
\end{align}
From inequalities (\ref{cond1}) and (\ref{nas2}) it follows that it now suffices to find an $x_*$ such that for all $x>x_*$ the following inequality holds: 
\begin{align}\label{ghj}
\frac{1}{\nu+1/2+x}<\frac{1+\gamma(2\nu+1)}{2(\nu+1)(1+\gamma x)}.   
\end{align}
On rearranging inequality (\ref{ghj}), we conclude that 
$u_{\nu,\gamma}'(x)>0$ for all $x>x_*$, where
\begin{equation}\label{xstar}
x_*=\frac{1}{2}(2\nu+1)^2+\frac{(1-2\nu)(\nu+1)}{1-\gamma}.
\end{equation}
Note that $x_*<0$ for some values of $\nu>-1/2$ and $0\leq\gamma<1$; for such values of $\nu$ and $\gamma$ we have $u_{\nu,\gamma}'(x)>0$ for all $x>0$, meaning that inequality (\ref{mainbd}) is satisfied for all $x>0$, since $u_{\nu,\gamma}(0^+)>0$ for all $\nu>-1/2$, $0\leq\gamma<1$ (which is easily seen from an application of the limiting form (\ref{Itend0})).
Similarly, due to the bound $I_{\nu+1}(x)/I_\nu(x)<x/(2\nu+2)$, (for $x>0$, $\nu>-1$) of \cite{nasell2} we can deduce from (\ref{cond1}) that there exists an $x_{**}$ such that $u_{\nu,\gamma}'(x)>0$ for all $x\in(0,x_{**})$, where
\begin{equation*}
x_{**}=2\nu+1.    
\end{equation*}
It is easily seen that $x_{*}\leq x_{**}$ for $\nu\geq1/2$, $0\leq\gamma<1$, and it therefore follows that inequality (\ref{mainbd}) holds for all $x>0$ when $\nu\geq1/2$, $0\leq \gamma<1$. 

\vspace{2mm}

\noindent 2.  It now suffices to prove inequality (\ref{mainbd}) for $0<x\leq x_*$, with $x_*$ defined as in (\ref{xstar}), in the case $-1/2<\nu<1/2$, $0\leq \gamma<1$; note that $x_*>0$ in this case. To this end, we now derive a bound for the integral (\ref{besint}) which is convenient for the purposes of bounding this integral for $0<x\leq x_*$. Writing $\mathrm{e}^{-\gamma t}=\mathrm{e}^{(1-\gamma)t}\mathrm{e}^{-t}$ and then performing an integration by parts with the aid of the indefinite integral formula (\ref{intfor}) yields
\begin{align*}
\int_0^x \mathrm{e}^{-\gamma t}t^{\nu}I_{\nu}(t)\,\mathrm{d}t=\frac{\mathrm{e}^{-\gamma x}x^{\nu+1}}{2\nu+1}\big(I_\nu(x)+I_{\nu+1}(x)\big)-\frac{1-\gamma}{2\nu+1}\int_0^x\mathrm{e}^{-\gamma t}t^{\nu+1} \big(I_\nu(t)+I_{\nu+1}(t)\big)\,\mathrm{d}t,
\end{align*}
where we used that $\lim_{x\rightarrow0}x^{\nu+1}I_{\nu}(x)=0$ and $\lim_{x\rightarrow0}x^{\nu+1}I_{\nu+1}(x)=0$ for $\nu>-1/2$ (see (\ref{Itend0})). Applying the integral inequalities (\ref{lower1}) and (\ref{lower3}) now leads to the bound
\begin{align}
\int_0^x \mathrm{e}^{-\gamma t}t^{\nu}I_{\nu}(t)\,\mathrm{d}t&<\frac{\mathrm{e}^{-\gamma x}x^{\nu+1}}{2\nu+1}\bigg(I_\nu(x)+I_{\nu+1}(x)\nonumber\\
&\quad-(1-\gamma)\sum_{k=0}^\infty\gamma^k\big(I_{\nu+k+1}(x)+I_{\nu+k+2}(x)\big)\bigg). \label{need}
\end{align}
We now truncate the infinite series in (\ref{need}) to obtain a simpler bound. Retaining further terms does not seem to be helpful in terms of improving on the constant factor $2(\nu+1)+c_\nu$; see Remark \ref{remcom} for further comments. 
Our simpler bound is
\begin{align}
\int_0^x \mathrm{e}^{-\gamma t}t^{\nu}I_{\nu}(t)\,\mathrm{d}t&<\frac{\mathrm{e}^{-\gamma x}x^{\nu+1}}{2\nu+1}\Big(I_\nu(x)+I_{\nu+1}(x)-(1-\gamma)\big(I_{\nu+1}(x)+I_{\nu+2}(x)+\gamma I_{\nu+2}(x)\big)\Big)\nonumber\\
&=\frac{\mathrm{e}^{-\gamma x}x^{\nu+1}}{2\nu+1}\bigg(\bigg(\frac{2(\nu+1)}{x}+\gamma\bigg)I_{\nu+1}(x)+\gamma^2I_{\nu+2}(x)\bigg),\label{need2}
\end{align}
where we obtained the equality using the identity (\ref{Iidentity}).

\vspace{2mm}

\noindent 3. We now consider the function
\begin{align}\label{vfor}
v_{\nu,\gamma}(x)=(1-\gamma)\bigg(2(\nu+1)+\gamma x+\gamma^2\frac{xI_{\nu+2}(x)}{I_{\nu+1}(x)}\bigg).  
\end{align}
To complete the proof of inequality (\ref{mainbd}), we note that it follows from inequality (\ref{need2}) that it will suffice to prove that $v_{\nu,\gamma}(x)\leq 2\nu+2+c_\nu$ for $x\in(0,x_*]$. 
To establish this bound, we will exploit the fact that $v_{\nu,\gamma}(x)$ is an increasing function of $x$ on the interval $(0,x_*]$, which follows because, for $\nu>-1/2$, the ratio $I_{\nu+2}(x)/I_{\nu+1}(x)$ is an increasing function of $x$ on the interval $(0,\infty)$ (see \cite[inequality (15)]{amos}). Therefore, for $x\in(0,x_*]$,
\begin{align*}
v_{\nu,\gamma}(x)\leq v_{\nu,\gamma}(x_*)&<(1-\gamma)\bigg(2(\nu+1)+(\gamma+\gamma^2)\bigg(\frac{1}{2}(2\nu+1)^2+\frac{(1-2\nu)(\nu+1)}{1-\gamma}\bigg)\bigg) \\
&=2(\nu+1)-(\nu+1/2)\gamma+(\nu+1)(1-2\nu)\gamma^2-2(\nu+1/2)^2\gamma^3=:P_\nu(\gamma),
\end{align*}
where we bounded $I_{\nu+2}(x_*)/I_{\nu+1}(x_*)<1$ using inequality (\ref{Imon}).

The stationary points of the cubic $P_\nu(\gamma)$ are given by
\begin{equation}\label{stat}
\gamma_{\pm}=\frac{2(\nu+1)(1-2\nu)\pm\sqrt{4(\nu+1)^2(1-2\nu)^2-3(2\nu+1)^3}}{3(2\nu+1)^2},    
\end{equation}
and it is easily seen that $\gamma_+\geq1$ for $-1/2<\nu\leq0$. Moreover, since $P_\nu'(0^+)<0$, the negative root $\gamma_-$ must be a minimum. Thus, for $-1/2<\nu\leq0$, we have the bound \begin{align*}
P_\nu(\gamma)&\leq \max\{P_\nu(0),P_\nu(1)\}\\
&=\max\{2(\nu+1),\,2(1-2\nu)(v+1)\}\\
&\leq 2(1-2\nu)(v+1)=2(\nu+1)-4\nu(\nu+1),  
\end{align*}
which proves that inequality (\ref{mainbd}) holds for all $x>0$ and $0\leq\gamma<1$ in the case $-1/2<\nu\leq0$.

We now consider the case $0<\nu<1/2$. From (\ref{stat}) we see that there are no stationary points if $\Delta(\nu)=4(\nu+1)^2(1-2\nu)^2-3(2\nu+1)^3<0$. This quartic inequality can be solved numerically to reveal that $\Delta(\nu)<0$ for $0.036049<\nu<2.1966$ (here and in the  remainder of the proof, we report roots of polynomial equations to 5 significant figures). Thus, for $0.036049<\nu<1/2$ the cubic $P_\nu(\gamma)$ is a decreasing function of $\gamma$ (since $P_\nu'(0^+)<0$), and therefore $P_\nu(\gamma)\leq P_\nu(0)=2(\nu+1)$. This proves that inequality (\ref{mainbd}) holds for all $x>0$ and $0\leq\gamma<1$ when $0.036049<\nu<1/2$.

It remains to treat the case $0<\nu\leq0.036049$. We consider the function $Q_\nu(\gamma)=P_\nu(\gamma)-2(\nu+1)$, and prove that $Q_\nu(\gamma)\leq 0$ for all $0\leq\gamma<1$ when $0<\nu\leq0.036049$. We have $Q_\nu(0)=0$ and $Q_\nu(1)=-4\nu(\nu+1)<0$ where the inequality holds because $\nu>0$. It therefore remains to prove that $Q_\nu(\gamma+)\leq0$, where $\gamma_+$ is the local maximum of the cubic $Q_\nu(\gamma)$. Since $\gamma_+>0$, it suffices to prove that $f_{\gamma_+}(\nu)\leq0$, where
\begin{align*}
f_\gamma(\nu)=-(\nu+1/2)+(\nu+1)(1-2\nu)\gamma-2(\nu+1/2)^2\gamma^2.
\end{align*}
A tedious calculation shows that
\begin{align*}
f_{\gamma_+}(\nu)=-\frac{g(\nu)}{9(2\nu+1)^2},
\end{align*}
where
\begin{equation*}
g(\nu)=1+22\nu+42\nu^2+16\nu^3-8\nu^4-(\nu+1)(1-2\nu)\sqrt{4(\nu+1)^2(1-2\nu)^2-3(2\nu+1)^3}.  
\end{equation*}
A simple calculation reveals that the roots of the equation $g(\nu)=0$ are given as solutions to the equation
$\nu(2\nu+1)^3(4\nu^3-4\nu^2-15\nu-8)=0$. The real-valued solutions to this equation are given by $\nu=-1/2$, $\nu=0$ and $\nu=2.6787$. Moreover, a simple asymptotic analysis shows that $f_{\gamma_+}(\nu)\sim -4\nu$ as $\nu\rightarrow0$, so that $f_{\gamma_+}(0)<0$ and $f_{\gamma_+}'(0^+)<0$ for $0<\nu\leq0.036049$. We therefore conclude that $f_{\gamma_+}(\nu)\leq0$, meaning that $P_\nu(\gamma)\leq P_\nu(0)=2(\nu+1)$ for all $0\leq \gamma<1$ when $0<\nu\leq0.036049$. 
This completes the proof of inequality (\ref{mainbd}) for all $x>0$, $\nu>-1/2$ and $0\leq\gamma<1$.
\qed

\begin{remark}\label{remcom} If we were to retain all terms in the infinite series from inequality (\ref{need}), then the role of the function $v_{\nu,\gamma}(x)$, as given in equation (\ref{vfor}), would instead be taken by the function
\begin{equation*}
\tilde{v}_{\nu,\gamma}(x)= (1-\gamma)\bigg(\frac{xI_\nu(x)}{I_{\nu+1}(x)}+x-(1-\gamma)\sum_{k=0}^\infty\gamma^k\bigg(\frac{xI_{\nu+k+1}(x)}{I_{\nu+1}(x)}+\frac{xI_{\nu+k+2}(x)}{I_{\nu+1}(x)}\bigg)\bigg).  
\end{equation*}
Given that $\tilde{v}_{\nu,\gamma}(x)<v_{\nu,\gamma}(x)$ for all $x>0$ (when $\gamma>0$), one may hope to be able to get a tighter bound for $\tilde{v}_{\nu,\gamma}(x)$ than the bound $v_{\nu,\gamma}(x)\leq 2\nu+2+c_\nu$, for $-1/2<\nu<1/2$, $0\leq\gamma<1$, that is given in the proof of Theorem \ref{thm1}. However, it does not appear that this approach will lead to an improved bound that holds for all $-1/2<\nu<1/2$, $0\leq\gamma<1$.
To see why, we introduce the  functions
\begin{equation*}
S_{\nu,j}(\gamma)=(1-\gamma)^2\sum_{k=0}^\infty\gamma^kR_{\nu+k+j}\bigg(\frac{1}{2}(2\nu+1)^2+\frac{(1-2\nu)(\nu+1)}{1-\gamma}\bigg), \quad j=1,2,   
\end{equation*}
where $R_{\nu+k+j}(x)=xI_{\nu+k+j}(x)/I_{\nu+1}(x)$. Numerical experiments with \emph{Mathematica} strongly suggest that $S_{\nu,j}(\gamma)\rightarrow0$ as $\gamma\rightarrow1$, for fixed $\nu>-1/2$ and $j=1,2$. Assuming that this is the case, we would then have that $\tilde{v}_{\nu,\gamma}(x_*)\rightarrow2(1-2\nu)(\nu+1)=2\nu+2+c_\nu$ as $\gamma\rightarrow1$, for fixed $\nu>-1/2$ and $j=1,2$, meaning that no further improvement is gained from retaining further terms in the infinite series.
\end{remark}


We now prove Theorem \ref{thm2}. Our proof we will make use of the following recent refinement of inequality (\ref{Imon}) due to \cite[Proposition 1]{yt23}. Let $\nu>-1/2$ and $0<a_\nu<1$. Then the inequality
\begin{equation}\label{ytineq}
I_{\nu+1}(x)<(1-a_\nu)I_\nu(x)+a_\nu I_{\nu+2}(x)   
\end{equation}
holds for all $x>0$ if and only $0<a_\nu\leq(\nu+1/2)/(2\nu+2)$.

\vspace{3mm}

\noindent{\emph{Proof of Theorem \ref{thm2}.}} (i) 
We begin by noting that inequality (\ref{new1}) was proved in the case $\gamma=0$ by \cite[inequality (2.17)]{gaunt ineq3}, in which case inequality (\ref{new1}) holds without the need for the restriction $\nu\geq1/2$ .
We may therefore restrict our attention to the $0<\gamma<1$ case. We also remark that the condition $\nu>-(n+1)/2$ is required to ensure that all integrals given in this proof exist (see (\ref{Itend0})). With the aid of the differentiation formula (\ref{diffone}) we obtain that
\begin{align}
\frac{\mathrm{d}}{\mathrm{d}t}\big(\mathrm{e}^{-\gamma t}t^{\nu}I_{\nu+n+1}(t)\big)
&=\frac{\mathrm{d}}{\mathrm{d}t}\big(\mathrm{e}^{-\gamma t}t^{-(n+1)}\cdot t^{\nu+n+1}I_{\nu+n+1}(t)\big)   \nonumber\\
&=\mathrm{e}^{-\gamma t}t^{\nu}\big(-\gamma I_{\nu+n+1}(t)-(n+1)t^{-1}I_{\nu+n+1}(t)+I_{\nu+n}(t)\big)\nonumber\\
&=\mathrm{e}^{-\gamma t}t^\nu\bigg(\frac{2\nu+n+1}{2(\nu+n+1)}I_{\nu+n}(t)+\frac{n+1}{2(\nu+n+1)}I_{\nu+n+2}(t)\nonumber\\
&\quad-\gamma I_{\nu+n+1}(t)\bigg), \label{eqnbad}
\end{align}
where we used identity (\ref{Iidentity}) to obtain the final equality. By inequality (\ref{ytineq}) we have that, for all $t>0$,
\begin{align}\label{ineqbad}
I_{\nu+n+1}(t)<\frac{2\nu+n+1}{2(\nu+n+1)}I_{\nu+n}(t)+\frac{n+1}{2(\nu+n+1)}I_{\nu+n+2}(t),
\end{align}
provided $0<(n+1)/(2\nu+2n+2)\leq(\nu+n+1/2)/(2\nu+2n+2)$, which holds for $n>-1$ and $\nu\geq1/2$. Applying inequality (\ref{ineqbad}) to (\ref{eqnbad}) and then integrating both sides of the resulting inequality over $(0, x)$ and rearranging yields 
\begin{align}
\int_0^x \mathrm{e}^{-\gamma t}t^{\nu}I_{\nu+n}(t)\,\mathrm{d}t&< \frac{2(\nu +n+1)}{(2\nu +n+1)(1-\gamma)}\mathrm{e}^{-\gamma x} x^{\nu} I_{\nu +n+1} (x)\nonumber\\
&\quad-\frac{n+1}{(2\nu +n+1)(1-\gamma)}\int_0^x \mathrm{e}^{-\gamma t}t^{\nu}I_{\nu+n+2}(t)\,\mathrm{d}t \label{badca}
\end{align}
(the inequality here is strict because $0<\gamma<1$). From inequality (2.16) of \cite{gaunt ineq3} we have that
\begin{equation}\label{day}
\int_0^x \mathrm{e}^{-\gamma t}t^{\nu}I_{\nu+n+2}(t)\,\mathrm{d}t> \mathrm{e}^{-\gamma x}x^{\nu}I_{\nu+n+3}(x), \quad x>0,\, \nu>-\tfrac{1}{2}(n+3),\,n>-3. 
\end{equation}
On applying this  inequality to the integral on the right-hand side of (\ref{badca}) we now obtain inequality (\ref{new1}).

Now suppose that $\gamma=0$. If we further set $n=-1$, then by the differentiation formula (\ref{diffone}) and the limiting form (\ref{Itend0}) we have the equality $\int_0^x t^\nu I_{\nu-1}(t)\,\mathrm{d}t=x^\nu I_{\nu}(x)$, confirming the assertion that we have equality when $\gamma=0$ and $n=-1$ and $\nu>0$. 

Suppose now that $\gamma=0$ and $n<-1$. Then on setting $\gamma=0$ in (\ref{eqnbad}) (since $n<-1$ we now require that $\nu>-(n+1)$, so that $\nu+n+1>0$ to ensure that all terms in (\ref{eqnbad}) are well-defined) and then integrating over $(0,x)$ and rearranging we get that 
\begin{equation}\label{badbad}\int_0^x t^{\nu} I_{\nu +n} (t)\,\mathrm{d}t = \frac{2(\nu +n+1)}{2\nu +n+1} x^{\nu} I_{\nu +n+1} (x) - \frac{n+1}{2\nu +n+1} \int_0^x t^{\nu} I_{\nu +n +2} (t)\,\mathrm{d}t.
\end{equation}
Under the assumptions $n<-1$ and $\nu>-(n+1)/2$, the factor $(n+1)/(2\nu+n+1)$ is now negative, and so on applying inequality (\ref{day}) (for which we require that $n>-3$ and $\nu>-(n+3)/2$) to the integral on the right-hand side of (\ref{badbad}) we now get a lower bound for $\int_0^x t^{\nu} I_{\nu +n} (t)\,\mathrm{d}t$.

\vspace{2mm}

\noindent (ii) Since $\mathrm{e}^{-\gamma t}$ is a non-increasing function of $t$ (as $0\leq\gamma<1$) we have that
\begin{align}&\int_0^x\mathrm{e}^{-\gamma t}t^\nu I_{\nu+n}(t)\,\mathrm{d}t\geq\mathrm{e}^{-\gamma x}\int_0^xt^\nu I_{\nu+n}(t)\,\mathrm{d}t\nonumber\\
&\quad=\frac{\mathrm{e}^{-\gamma x}}{2\nu+n+1}\bigg(2(\nu +n+1) x^{\nu} I_{\nu +n+1} (x) - (n+1) \int_0^x t^{\nu} I_{\nu +n +2} (t)\,\mathrm{d}t\bigg),\label{bad3}
\end{align}
where we used equation (\ref{badbad}) to obtain the equality. Applying the upper bound (\ref{new1}) (with $\gamma=0$) to bound the integral on the right-hand side of (\ref{bad3}) now yields the lower bound (\ref{lower4}).

\vspace{2mm}

\noindent (iii) By the limiting form (\ref{Itend0}) we have that, as $x\rightarrow0$, 
\begin{align*}
\int_0^x \mathrm{e}^{-\gamma t}t^{\nu}I_{\nu}(t)\,\mathrm{d}t&\sim\frac{x^{2\nu+1}}{2^\nu(2\nu+1)\Gamma(\nu+1)},  \\  
\frac{2(\nu+1)}{(2\nu+1)(1-\gamma)}\mathrm{e}^{-\gamma x}x^\nu I_{\nu+1}(x)&\sim\frac{x^{2\nu+1}}{2^\nu(2\nu+1)(1-\gamma)\Gamma(\nu+1)}.    
\end{align*}
The tightness of inequality (\ref{new1}) as $x\rightarrow0$ when $\gamma=0$, and the tightness of inequality (\ref{lower4}) as $x\rightarrow0$ now follows.
The tightness of inequality (\ref{new1}) as $x\rightarrow\infty$ and the tightness of inequality (\ref{lower4}) as $x\rightarrow\infty$ when $\gamma=0$ follows from the same argument as used in part (iv) of the proof of Theorem \ref{cor1}. \qed

\vspace{3mm}

\noindent{\emph{Proof of Proposition \ref{prop1}.}} (i) Since $\mu\geq\nu\geq1/2$, we have that, for $x>0$,
\begin{align*}
\int_0^x \mathrm{e}^{-\gamma t}t^{\mu}I_{\nu}(t)\,\mathrm{d}t\leq x^{\mu-\nu}\int_0^x \mathrm{e}^{-\gamma t}t^{\nu}I_{\nu}(t)\,\mathrm{d}t<\frac{1}{1-\gamma}\mathrm{e}^{-\gamma x}x^\mu I_\nu(x), 
\end{align*}
where we used inequality  (\ref{gau2}) to obtain the second inequality.

\vspace{2mm}

\noindent (ii) By a simple asymptotic analysis using (\ref{Itendinfinity}) and integration by parts we have that, for $\mu+\nu>-1$ and $0\leq\gamma<1$, 
\begin{equation}\label{eqeq12} \int_0^x \mathrm{e}^{-\gamma t}t^\mu  I_{\nu}(t)\,\mathrm{d}t\sim \frac{x^{\mu-1/2}\mathrm{e}^{(1-\gamma)x}}{\sqrt{2\pi}(1-\gamma)}\bigg\{1-\bigg(\frac{4\nu^2-1}{8}+\frac{\mu-1/2}{1-\gamma}\bigg)\frac{1}{x}\bigg)\bigg\}, \quad x\rightarrow\infty,
\end{equation}
in which we retained the second term in the expansion (\ref{Itendinfinity}) to get the second term in the asymptotic expansion of the integral in (\ref{eqeq12}). In comparison, it is immediate from (\ref{Itendinfinity}) that
\begin{align}\label{inff}
\frac{1}{1-\gamma}\mathrm{e}^{-\gamma x}x^\mu I_\nu(x)\sim \frac{x^{\mu-1/2}\mathrm{e}^{(1-\gamma)x}}{\sqrt{2\pi}(1-\gamma)}\bigg(1-\frac{4\nu^2-1}{8x}\bigg),\quad x\rightarrow\infty.  
\end{align}
It can be seen that the second term in the expansion (\ref{eqeq12}) is greater than the second term in expansion of (\ref{inff}) if $(\mu-1/2)/(1-\gamma)<0$, that is if $\mu<1/2$. Hence, when $\mu<1/2$ inequality (\ref{lar}) must hold for sufficiently large $x>0$. \qed

\vspace{3mm}

We now prove Corollary \ref{cor2}, for which we will require the following inequalities for products of modified Bessel functions.  Inequality (\ref{prod1}) is  given in the proof of Theorem 5 of \cite{gaunt ineq2}, and is a simple consequence of Theorem 4.1 of \cite{hartman}. Inequality (\ref{prod2}) is given in \cite[Lemma 3]{gaunt ineq2}. Further results and bounds for the product $I_{\nu}(x)K_{\nu}(x)$ can be found in \cite{baricz,baricz2}.  
For $x\geq0$, 
\begin{equation}\label{prod1}0\leq xK_{\nu}(x)I_\nu(x)<\frac{1}{2}, \quad \nu>\tfrac{1}{2},
\end{equation}
and
\begin{equation}\label{prod2}\frac{1}{2}< xK_{\nu+1}(x)I_\nu(x)\leq1, \quad \nu\geq-\tfrac{1}{2}.
\end{equation}

\vspace{3mm}

\noindent{\emph{Proof of Corollary \ref{cor2}.}} (i) An application of inequality (\ref{mainbd}) followed by inequality (\ref{prod2}) yields
\begin{align*}\frac{\mathrm{e}^{-\beta x}K_{\nu+2}(x)}{x^{\nu-1}}\int_0^x \mathrm{e}^{\beta t}t^{\nu}I_\nu(t)\,\mathrm{d}t&<\frac{2(\nu+1)+c_\nu}{(2\nu+1)(1+\beta)}x K_{\nu+2}(x) I_{\nu+1}(x)\leq\frac{2(\nu+1)+c_\nu}{(2\nu+1)(1+\beta)}.
\end{align*}
\noindent (ii) This time we may apply inequality (\ref{prod1}) to obtain
\begin{align*}\frac{\mathrm{e}^{-\beta x}K_{\nu+1}(x)}{x^{\nu-1}}\int_0^x \mathrm{e}^{\beta t}t^{\nu}I_\nu(t)\,\mathrm{d}t&<\frac{2(\nu+1)+c_\nu}{(2\nu+1)(1+\beta)}x K_{\nu+1}(x) I_{\nu+1}(x)<\frac{\nu+1+c_\nu/2}{(2\nu+1)(1+\beta)}.
\end{align*}
\noindent (iii) This follows from part (ii) and inequality (\ref{cake}). 
\qed

\appendix

\section{Basic properties of modified Bessel functions}\label{appa}

In this appendix, we state some basic properties of modified Bessel functions that are used in this paper.  All formulas can be found in the standard reference \cite{olver}, except for the inequalities.

The modified Bessel functions of the first kind $I_\nu(x)$ and second kind $K_\nu(x)$ are defined, for $\nu\in\mathbb{R}$ and $x>0$, by
\begin{equation*}I_{\nu} (x) =  \sum_{k=0}^{\infty} \frac{(x/2)^{\nu+2k}}{\Gamma(\nu +k+1) k!}, \quad K_\nu(x)=\int_0^\infty \mathrm{e}^{-x\cosh(t)}\cosh(\nu t)\,\mathrm{d}t.
\end{equation*}
For $x>0$, the functions $I_{\nu}(x)$ and $K_{\nu}(x)$ are positive for $\nu\geq-1$ and all $\nu\in\mathbb{R}$, respectively.  The modified Bessel function $I_\nu(x)$ satisfies the following identity and differentiation formula:
\begin{align}
\label{Iidentity}I_{\nu +1} (x) &= I_{\nu -1} (x) - \frac{2\nu}{x} I_{\nu} (x), \\
\label{diffone}\frac{\mathrm{d}}{\mathrm{d}x} (x^{\nu} I_{\nu} (x) ) &= x^{\nu} I_{\nu -1} (x).
\end{align}
The following indefinite integration formula also holds:
\begin{equation}\label{intfor}\int \mathrm{e}^{-x}x^\nu I_\nu(x)\,\mathrm{d}x=\frac{\mathrm{e}^{-x}x^{\nu+1}}{2\nu+1}\big(I_\nu(x)+I_{\nu+1}(x)\big), \quad x>0,\; \nu>-\tfrac{1}{2}.
\end{equation}
The modified Bessel function of the first kind possess the following asymptotic behaviour:
\begin{align}\label{Itend0}I_{\nu} (x) &\sim \frac{(\frac{1}{2}x)^\nu}{\Gamma(\nu +1)}\bigg(1+\frac{x^2}{4(\nu+1)}\bigg), 
\quad x \downarrow 0, \: \nu\notin\{-1,-2,-3,\ldots\}, \\
\label{Itendinfinity}I_{\nu} (x) &\sim \frac{\mathrm{e}^{x}}{\sqrt{2\pi x}}\bigg(1-\frac{4\nu^2-1}{8x}\bigg), \quad x \rightarrow\infty, \: \nu\in\mathbb{R}.
\end{align}
For $x > 0$, the following inequalities are satisfied:
\begin{align}\label{Imon}I_{\nu+1} (x) &< I_{\nu} (x), \quad \nu \geq -\tfrac{1}{2}, \\
\label{cake}K_{\nu+1} (x) &> K_{\nu} (x), \quad \nu > -\tfrac{1}{2}.
\end{align}
 Inequality (\ref{Imon}) is given in \cite{jones} and \cite{nasell}, which extended a result of \cite{soni}. Inequality (\ref{cake}) can be found in \cite{ifant}. A survey of inequalities for modified Bessel functions is given by  \cite{baricz24}. For refinements of inequalities (\ref{Imon}) and (\ref{cake}), see \cite{segura,segura2} and references therein. 

\section*{Acknowledgements}
The author is funded in part by EPSRC grant EP/Y008650/1 and EPSRC grant UKRI068.


\footnotesize

\begin{thebibliography}{99}
\addcontentsline{toc}{section}{References}

\bibitem{amos} Amos, D. E. Computation of modified Bessel functions and their ratios. \emph{Math. Comp.} $\mathbf{28}$ (1974), 239--251.

\bibitem{aaps17} Arras, B., Azmoodeh, E., Poly, G. and Swan, Y. A bound on the 2-Wasserstein distance between linear combinations of independent random variables. \emph{Stoch. Proc. Appl.} $\mathbf{129}$ (2019), 2341--2375. 

\bibitem{aet21} Azmoodeh, E., Eichelsbacher, P. and Th\"{a}le, C. Optimal Variance-Gamma approximation on the second Wiener chaos. \emph{J. Funct. Anal.} $\mathbf{282}$ (2022), Art.\ 109450.


\bibitem{ag18} Azmoodeh, E. and Gasbarra, D. On a new Sheffer class of polynomials related to normal product distribution. \emph{Theor. Probab. Math. Statist.} $\mathbf{98}$ (2019), 51--71. 

\bibitem{azmooden} Azmoodeh, E., Peccati, G. and Poly, G.  Convergence towards linear combinations of chi-squared random variables: a Malliavin-based approach.  \emph{S\'{e}minaire de Probabilit\'{e}s} XLVII (special volume in memory of Marc Yor) (2015), 339--367.

\bibitem{bt17} Bai, S. and Taqqu, M.  Behavior of the generalized Rosenblatt process at extreme critical exponent values. \emph{Ann. Probab.} $\mathbf{45}$ (2017), 1278--1324.


\bibitem{baricz} Baricz, \'{A}.  On a product of modified Bessel functions.  \emph{Proc. Amer. Math. Soc.} $\mathbf{137}$ (2009),  189--193.

\bibitem{baricz24} Baricz, \'{A}.  Bounds for modified Bessel functions of the first and second kinds. \emph{Proc. Edinburgh Math. Soc.} $\mathbf{53}$ (2010), 575--599.

\bibitem{baricz2} Baricz, \'{A}., Jankov Ma\v{s}irevi\'{c}, D. J., Ponnusamy, S. and Singh, S.  Bounds for the product of modified Bessel functions.  \emph{Aequat. Math.} $\mathbf{90}$ (2016), 859--870.

\bibitem{bs09} Baricz, \'{A}. and Sun, Y. New bounds for the generalized Marcum $Q$-function. \emph{IEEE Trans. Info. Th.} $\mathbf{55}$ (2009), 3091--3100.

\bibitem{baricz3} Baricz, \'{A}. and Sun, Y.  Bounds for the generalized Marcum $Q$-function. \emph{Appl. Math. Comput.} $\mathbf{217}$ (2010),  2238--2250.



\bibitem{chen} Chen, L. H. Y., Goldstein, L. and Shao, Q.--M.  \emph{Normal Approximation by Stein's Method.} Springer, 2011.


\bibitem{et15} Eichelsbacher, P. and Th\"{a}le, C.  Malliavin-Stein method for Variance-Gamma approximation on Wiener space.  	\emph{Electron. J. Probab.} $\mathbf{20}$ no. 123 (2015).

\bibitem{vg review} Fischer, A., Gaunt, R. E. and Sarantsev, A. The Variance-Gamma Distribution: A Review. To appear in \emph{Stat. Sci.}, 2025+.



\bibitem{gaunt vg} Gaunt, R. E.  Variance-Gamma approximation via Stein's method.  \emph{Electron. J. Probab.} $\mathbf{19}$ no. 38 (2014).

\bibitem{gaunt ineq1} Gaunt, R. E.  Inequalities for modified Bessel functions and their integrals.  \emph{J. Math. Anal. Appl.} $\mathbf{420}$ (2014),  373--386.

\bibitem{gaunt ineq2} Gaunt, R. E.  Uniform bounds for expressions involving modified Bessel functions.  \emph{Math. Inequal. Appl.} $\mathbf{19}$ (2016), 1003--1012.

\bibitem{gaunt ineq3} Gaunt, R. E. Inequalities for integrals of modified Bessel functions and expressions involving them.  \emph{J. Math. Anal. Appl.} $\mathbf{462}$ (2018),  172--190.


\bibitem{gaunt ineq6} Gaunt, R. E. Inequalities for some integrals involving modified Bessel functions. \emph{Proc. Amer. Math. Soc.} $\mathbf{147}$ (2019), 2937--2951.

\bibitem{gaunt vg2} Gaunt, R. E. Wasserstein and Kolmogorov error bounds for variance-gamma approximation via Stein's method I. \emph{J. Theor. Probab.} $\mathbf{33}$ (2020), 465--505.

\bibitem{gaunt2021} Gaunt, R. E. Bounds for an integral of the modified Bessel function of the first kind and expressions involving it. \emph{J. Math. Anal. Appl.} $\mathbf{502}$ (2021), Art.\ 125216.

\bibitem{gaunt vg3}  Gaunt, R. E. Stein factors for variance-gamma approximation in the Wasserstein and Kolmogorov distances. \emph{J. Math. Anal. Appl.} $\mathbf{514}$ (2022), Art.\ 126274.


\bibitem{hartman} Hartman, P. On the products of solutions of second order disconjugate differential equations and the
Whittaker differential equation. \emph{SIAM J. Math. Anal.} $\mathbf{8}$ (1977),  558--571.

\bibitem{ifant} Ifantis, E. K. and Siafarikas, P. D. Bounds for modified bessel functions. \emph{Rend. Circ. Mat. Palermo} $\mathbf{40}$ (1991),  347--356.

\bibitem{jones} Jones, A. L.  An extension of an inequality involving modified Bessel functions.  \emph{J. Math. Phys. Camb.} $\mathbf{47}$ (1968),  220--221.





\bibitem{nasell} N\r{a}sell, I. Inequalities for Modified Bessel Functions. \emph{Math. Comput.} $\mathbf{28}$ (1974),  253--256. 

\bibitem{nasell2} N\r{a}sell, I. Rational bounds for ratios of modified Bessel functions. \emph{SIAM J. Math. Anal.} $\mathbf{9}$ (1978), 1--11.

\bibitem{np12} Nourdin, I. and Peccati, G. Normal approximations with Malliavin calculus: from Stein's method to universality. Vol. 192. Cambridge University Press, 2012.

\bibitem{olver} Olver, F. W. J., Lozier, D. W., Boisvert, R. F. and Clark, C. W.  \emph{NIST Handbook of Mathematical Functions.} Cambridge University Press, 2010.




\bibitem{segura} Segura, J.  Bounds for ratios of modified Bessel functions and associated Tur\'{a}n-type inequalities.  \emph{J. Math. Anal. Appl.} $\mathbf{374}$ (2011), 516--528.

\bibitem{segura2} Segura, J. Simple bounds with best possible accuracy for ratios of modified Bessel functions. \emph{J. Math. Anal. Appl.} $\mathbf{502}$ (2023), Art.\ 127211.

\bibitem{soni} Soni, R. P.  On an inequality for modified Bessel functions. \emph{J. Math. Phys.} $\mathbf{44}$ (1965),  406--407.

\bibitem{stein} Stein, C.  A bound for the error in the normal approximation to the the distribution of a sum of dependent random variables.  In \emph{Proc. Sixth Berkeley Symp. Math. Statis. Prob.} (1972), vol. 2, Univ. California Press, Berkeley,  583--602.

\bibitem{yt23} Yang, Z.--H. and Tian, J.--F. Convexity of ratios of the modified Bessel functions of the first kind with applications. \emph{Rev. Mat. Complut.} $\mathbf{36}$ (2023), 799--825.

\end{thebibliography}
\end{document}